\pdfoutput=1

 \documentclass{amsart}

\usepackage{amssymb}

\usepackage[centertags]{amsmath}
\usepackage{amsfonts}
\usepackage{amssymb}
\usepackage{amsthm}
\usepackage{ subfigure}
\usepackage{graphicx}
\usepackage{fullpage }
\usepackage{color}

\newtheorem{theorem}{Theorem}

\newtheorem{lemma}[theorem]{Lemma}

\theoremstyle{definition}
\newtheorem{definition}[theorem]{Definition}

\theoremstyle{remark}
\newtheorem{remark}[theorem]{Remark}
\newtheorem{example}[theorem]{Example}

\newcommand{\F}{\mathcal{F}}
\newcommand{\calL}{\mathcal{L}}

\begin{document}



\title{Partial duality and Bollob\'as and Riordan's ribbon graph polynomial}


\author{Iain Moffatt}
\address{Department of Mathematics and Statistics,  University of South Alabama, Mobile, AL 36688, USA. }

\email{imoffatt@jaguar1.usouthal.edu}

\begin{abstract}
Recently S.~Chmutov introduced a generalization of the dual of a ribbon graph (or equivalently an embedded graph) and proved a relation between Bollob\'as and Riordan's ribbon graph polynomial of a ribbon graph and of its generalized duals. Here I show that the duality relation satisfied by the ribbon graph polynomial can be understood in terms of knot theory and I give a simple proof of the relation which used the homfly polynomial of a knot. 
\end{abstract}


\date{\today}

\maketitle









\section{Introduction and motivation}\label{s.intro}

Recently, there has been a lot of interest in connections between knots and ribbon graphs
(\cite{CP,CV,Ch1,Da,Mo1,Mo2}). In particular, there are  various constructions which realize the Jones polynomial of a link as an evaluation of Bollob\'as and Riordan's ribbon graph  polynomial (defined in \cite{BR1,BR}) of an associated signed ribbon graph. In \cite{CP}, Chmutov and Pak proved that the Jones polynomial of a virtual link or a link in a thickened surface is an evaluation of the signed ribbon graph  polynomial. In other work in this area,  Dasbach et. al. in  \cite{Da} showed how to construct a (non-signed) ribbon graph  from  from a (not necessarily alternating) link diagram with the property that the Jones polynomial is an evaluation of the ribbon graph 
 polynomial of the ribbon graph. 
 Given the similarity between these two results, as they both relate the Jones and ribbon graph polynomials,
  it is natural to look for a connection between them. This question was first answered in \cite{Mo2} where I defined an ``unsigning'' procedure  which took in a signed plane graph and gave out a non-signed ribbon graph. Chmutov has also considered the relationship between the ribbon graph models for the Jones polynomial, particularly between those in \cite{CP} and \cite{CV} (In \cite{CV} Chmutov and Voltz extended the results of Dasbach et. al. from \cite{Da} to virtual links).  In the process, Chmutov defined a generalized duality for ribbon graphs (which I call ``partial duality'' here\footnote{With thanks to Dan Archdeacon for suggesting the name ``partial duality''.}) of which my unsigning is a special case (as was observed by Chmutov in \cite{Ch1}). Chmutov not only showed that his partial duality connected ribbon graph models for the Jones polynomial, but also that it has desirable properties with respect to the signed ribbon graph  polynomial. These desirable properties generalize the well known behavior of the Tutte polynomial under duality.  In this paper I am interested in this partial duality and the ribbon graph  polynomial.

The partial dual $G^A$ of a ribbon graph $G$ is constructed by forming the dual of a ribbon graph only along the edges in  $A\subseteq E(G)$, as described in Subsection~\ref{ss.pd} below. Since $G^{E(G)}=G^*$, Poincar\'{e} duality is a special case of Chmutov's partial duality. Chmutov proved that, up to a normalization,  the signed ribbon graph  polynomials of $G$ and $G^A$ are  equal when $xyz^2=1$. In \cite{Mo1}, I used the fact that the homfly polynomial determines the ribbon graph  polynomial to prove that the ribbon graph  polynomials of $G$ and $G^*$ are equal, again up to a normalization, and again  along the surface $xyz^2=1$. Connecting the facts that Chmutov's duality relation holds along $xyz^2=1$; the homfly polynomial determines the ribbon graph  polynomial along $xyz^2$; and a special case ($A=E(G)$) of Chmutov's duality theorem has a simple proof through knot theory, one naturally suspects that  Chmutov's duality theorem can be understood in terms of knot theory. Here I show that this is indeed the case, and provide a proof of Chmutov's duality relation using knot theory. Showing that there is a knot theoretical foundation for this result offers a new understanding of the underlying structures of duality and the ribbon graph polynomial. 

The argument I use to  prove Chmutov's duality theorem is essentially: Step 1: the homfly determines the signed ribbon graph  polynomial; Step 2: the links associated with $G$ and $G^A$ have the same homfly polynomial.  This two step argument is structured in this paper in the following way.
 In Section~\ref{s.pd} I define the the partial dual of a signed, orientable ribbon graph. In Section~\ref{s.poly}, I review the definitions of the signed ribbon graph  polynomial and the homfly polynomial. I then go on to show how the homfly polynomial determines the signed ribbon graph  polynomial along $xyz^2=1$. This is split between sections in Section~\ref{ss.hbr}, where I review results from \cite{Mo1}, and  Section~\ref{ss.kt}, where I express the signed ribbon graph  polynomial in terms of the homfly polynomial and reformulate Chmutov's duality theorem. Finally, in Section~\ref{ss.proof}, I give a simple proof of the knot theoretic reformulation of the duality theorem.  

I would like to thank Tom Zaslavsky for encouraging me to write down these results.

\section{The partial dual of a ribbon graph}\label{s.pd}

\subsection{Ribbon graphs}

 Roughly speaking, a ribbon graph is a `topological graph' formed by using disks as vertices and ribbons $I\times I$ as edges. Ribbon graphs provide a convenient description of cellularly embedded graphs (a cellularly embedded graph  is an  embedded graph with the property that each of its faces is a 2-cell).

\begin{definition}
 A {\em ribbon graph} $G=(V(G),E(G))$ is
an surface with boundary represented as the union of
 closed  disks (called {\em vertices}) and   ribbons $I \times I$, where $I=[0,1]$ is the unit interval, (called {\em edges})  such that
\begin{enumerate}
\item the vertices and edges intersect in disjoint  line segments
$\{0,1\} \times I$;
\item each such line segment lies on the boundary of precisely one
vertex and precisely one edge;
\item every edge contains exactly two such line segments.
\end{enumerate}
A ribbon graph is said to be {\em orientable} if its underlying surface is orientable.

\end{definition}

A ribbon graph $G$ is said to be {\em signed} if it is equipped with a  a mapping from its edge set $E(G)$ to $ \{+,-\}$ (so a sign $+$ or $-$ is assigned to each edge of $G$). 

Ribbon graphs are considered up to homeomorphisms of the surface that preserve the vertex-edge structure. Some signed ribbon graphs are shown in examples~\ref{ex.dex1} and \ref{ex.dex2}.

 It is often convenient to label the edges of ribbon graphs. I will often abuse notation and identify an edge with its  unique label. At times I will also abuse notation and use $e$ to denote an edge of a ribbon graph and the label of that edge.

It is well known that ribbon graphs are equivalent to cellularly embedded graphs (considered up to  homeomorphism of the surface).  Details of equivalence of ribbon graphs and cellularly embedded graphs can be found in \cite{GT}, for example. Here I will work primarily in the language of ribbon graphs, rather than embedded graphs, as the topology of ribbon graphs is particularly convenient for my purposes. 

In this paper I will be primarily interested in orientable ribbon graphs.  An orientable ribbon graph is equivalent to a graph cellularly embedded in an orientable surface and is also 
equivalent to a  combinatorial map (that is a graph equipped with a cyclic order of the incident half-edges at each vertex). The restriction here  to orientable ribbon graphs is due to that fact that, at the time of writing, the homfly polynomial of a link in a thickened non-orientable surface has yet to be defined. It should be emphasized that all of the graph theoretical constructions used  in this paper do work for non-orientable ribbon graphs. Also, I expect that the knot theoretic methods used in this paper would extend to the non-orientable case with a suitable definition of the homfly polynomial of a link in a thickened non-orientable surface.

\subsection{Arrow presentations}

In order to define partial duality it will convenient to describe ribbon graphs using arrow presentations. Arrow presentations provide a useful combinatorial description of a ribbon graph. 
 \begin{definition}
From \cite{Ch1}, an {\em arrow presentation} consists of  a set of circles, called {\em cycles},  equipped  with a set of disjoint,  labelled arrows marked along their perimeters. Each label appears on precisely two arrows.  Two arrow presentations are considered equivalent if one can be obtained from the other by reversing the direction of all of the marking arrows which belong to some subset of labels or by changing the set of labels used.

An arrow presentation is said to be {\em signed} if there is a mapping  from the set of labels of the arrows to $\{+.-\}$.
\end{definition}

Arrow presentations and ribbon graphs are known to be equivalent (see \cite{GT}). This equivalence also holds for signed arrow presentations and signed ribbon graphs.  I will now describe how to move  between equivalent arrow presentations and ribbon graphs.

\begin{figure}
\[\includegraphics[height=15mm]{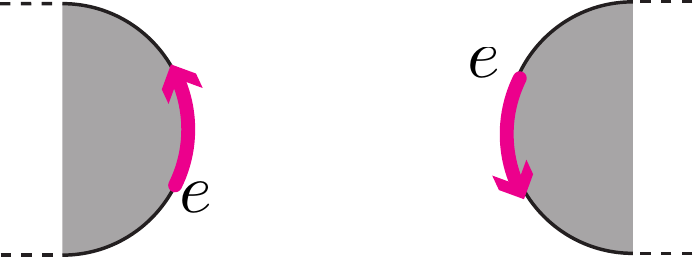}  
\raisebox{6mm}{\hspace{5mm}\includegraphics[width=11mm]{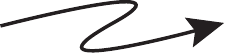}}
\includegraphics[height=15mm]{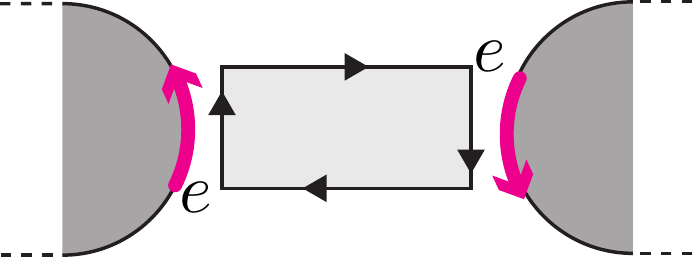}  
\raisebox{6mm}{\hspace{5mm}\includegraphics[width=11mm]{arrow}}
\includegraphics[height=15mm]{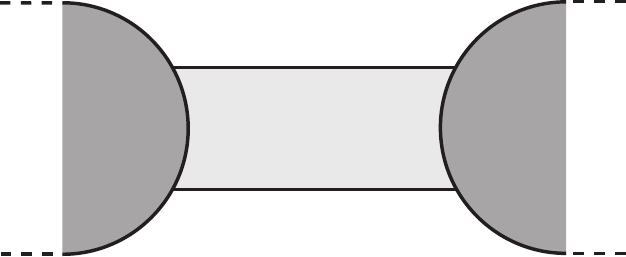} \]
\caption{Constructing a ribbon graph from an arrow presentation. }
\label{arrows}
\end{figure}

A ribbon graph can be obtained from an arrow presentation by viewing each cycle of the arrow presentation as the boundary of a disk that becomes a vertex of the ribbon graph. 
 Edges (which are 2-cells $I\times I$) are then added to the vertices (which are disks) by taking one edge $I\times I$ for each distinct label of the marking arrows then orienting the boundaries of the edges arbitrarily. Each edge  is then attached to one or two vertices by identifying each of the arcs $\{0\}\times I$ and $\{1\}\times I$ on the boundary of the edge with two arrows that have the same label. The edges are attached so that the orientation on the boundary of an edge agrees with the direction of the arrow. Moreover, exactly one arc on one edge is attached to each arrow. 
The process of attaching an edge is shown graphically in  Figure~\ref{arrows}.

Conversely, every ribbon graph  gives rise to an arrow presentation.  To describe a ribbon graph $G$ as an arrow presentation, start by arbitrarily  labelling and orienting the boundary of  each edge of $G$.  On the arcs $\{0\}\times I$ and $\{1\}\times I$, where an edge intersects a vertex, place a marked arrow on the vertex disk, labelling the arrow with the label of the edge it meets and directing the arrow consistently with the orientation of the boundary of the edge. The boundaries of the vertex set marked with these labelled arrows give the arrow marked cycles of an arrow presentation. 

If the arrow presentation is signed, then the edges of the corresponding ribbon graph naturally inherit signs from the labels of the arrows that the edges were attached to. Conversely, if a ribbon graph is signed, then the corresponding arrow presentation naturally inherits signs by associating the sign of each edge with the labels of the arrows it gives rise to. Thus, signed ribbon graphs are equivalent to signed arrow presentations. 

\begin{example} This is an  example of the equivalence between signed arrow presentations and signed ribbon graphs.  The labels $1$, $2$ and $3$ are arbitrary. Note that the ribbon graph is non-orientable.
\[ \raisebox{4mm}{\includegraphics[width=35mm]{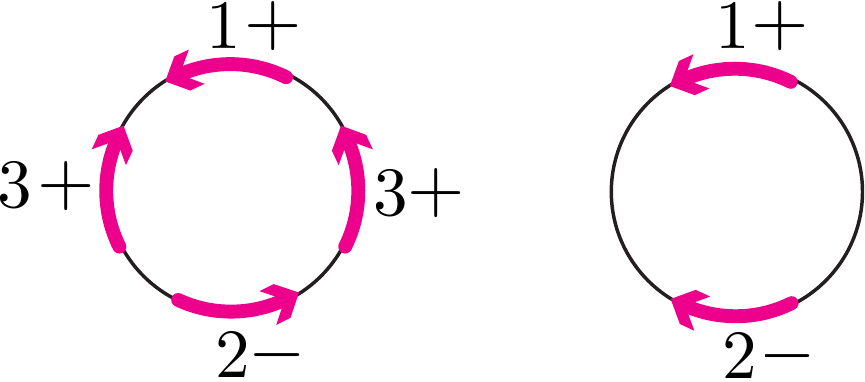}}\quad
\raisebox{10mm}{\includegraphics[width=18mm]{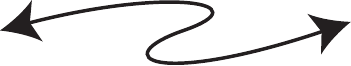}} \quad
 \includegraphics[width=40mm]{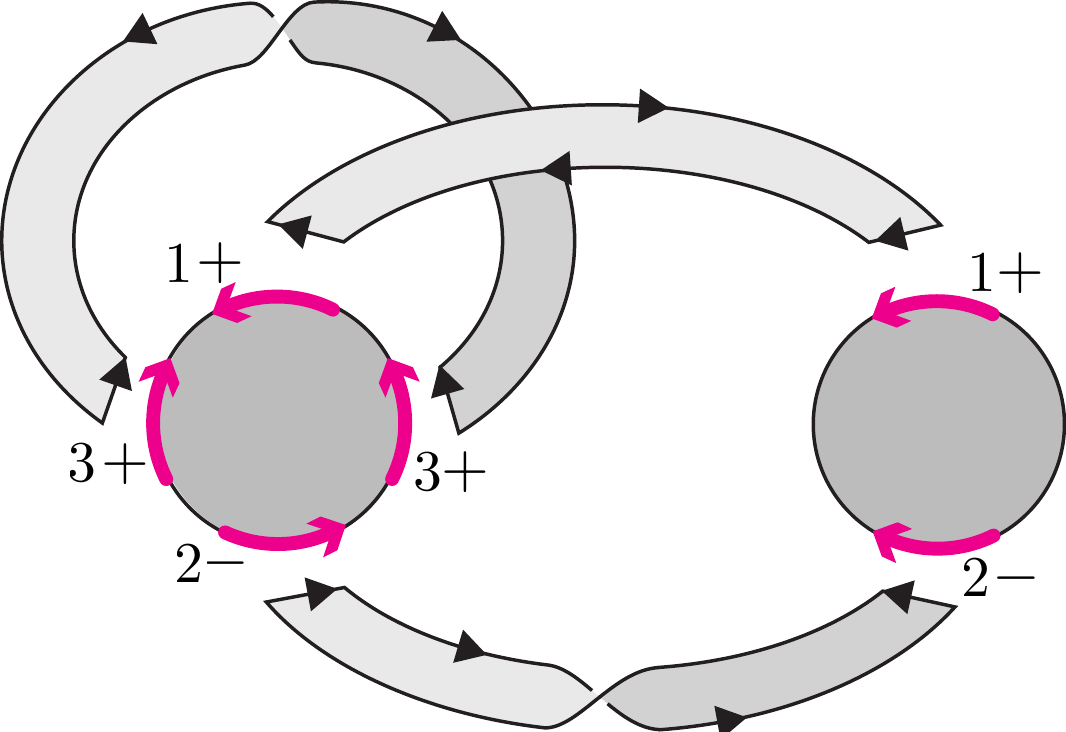} \quad
 \raisebox{10mm}{\includegraphics[width=18mm]{doublearrow}}\quad
\raisebox{1mm}{ \includegraphics[width=40mm]{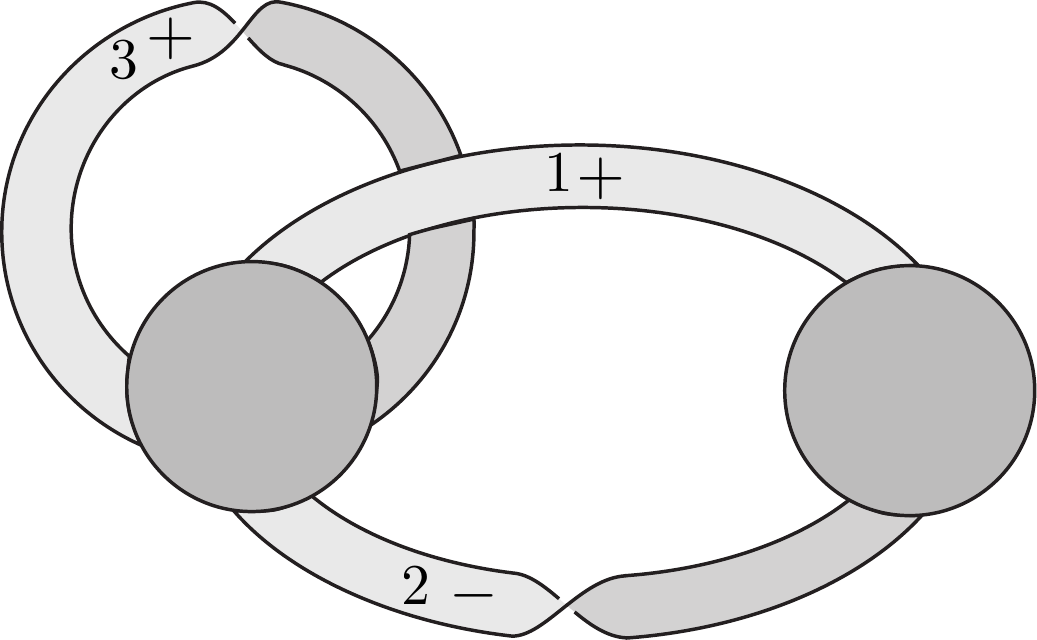}}.\]
\end{example}

\subsection{Partial duality}\label{ss.pd}

I will now give the definition of Chmutov's partial duality. The procedure used in the definition below starts with a signed ribbon graph and  a subset of edges and constructs a signed arrow presentation. The signed ribbon graph corresponding to the signed arrow presentation is a partial dual of the original signed ribbon graph.  

\begin{definition}
Let $G$ be a signed ribbon graph and $A\subseteq E(G)$. Arbitrarily orient and label each of the edges of $G$ (the orientation need not extend to an orientation of the ribbon graph). The boundary components of the spanning ribbon sub-graph $(V(G), A)$ of $G$ meet the edges of $G$ in disjoint arcs (where the spanning ribbon sub-graph is naturally embedded in $G$). On each of these arcs, place an arrow which points in the direction of the orientation of the edge  and is labelled by the edge it meets. Associate a sign to each label in the following way: if $e$ is a label of an edge of $G$ with sign $\varepsilon$, then the arrow labelled by $e$ has sign $-\varepsilon$ if the edge is in $A$, and has sign $\varepsilon$ otherwise. 
The resulting decorated boundary components of the spanning ribbon sub-graph $(V(G), A)$ define an signed arrow presentation. 
The signed ribbon graph corresponding to this signed arrow presentation is the {\em partial dual} $G^A$ of $G$.
\end{definition}

\begin{example}\label{ex.dex1} 
The signed ribbon graph $G$ equipped with and arbitrary labelling and orientation of its edges is shown in Step 1. In this example $A=\{2,3\}$. The marked spanning ribbon sub-graph $(V(G),A)$ is shown in Step 2 (note the change of signs of the edges in $A$). The boundary components of this give a signed arrow presentation, shown in Step 3. The corresponding signed ribbon graph is shown in Step 4. This is the partial dual $G^{\{2,3\}}$ of $G$.

 \begin{center}
\begin{tabular}{ccc}
 \includegraphics[width=5cm]{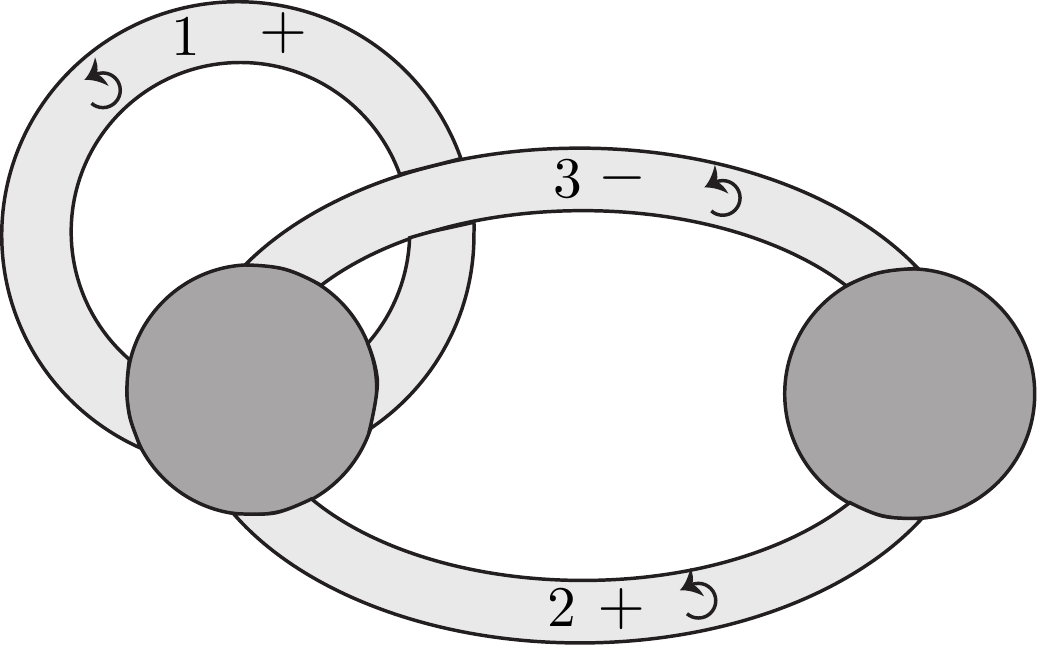}  & \hspace{2cm} & \includegraphics[width=5cm]{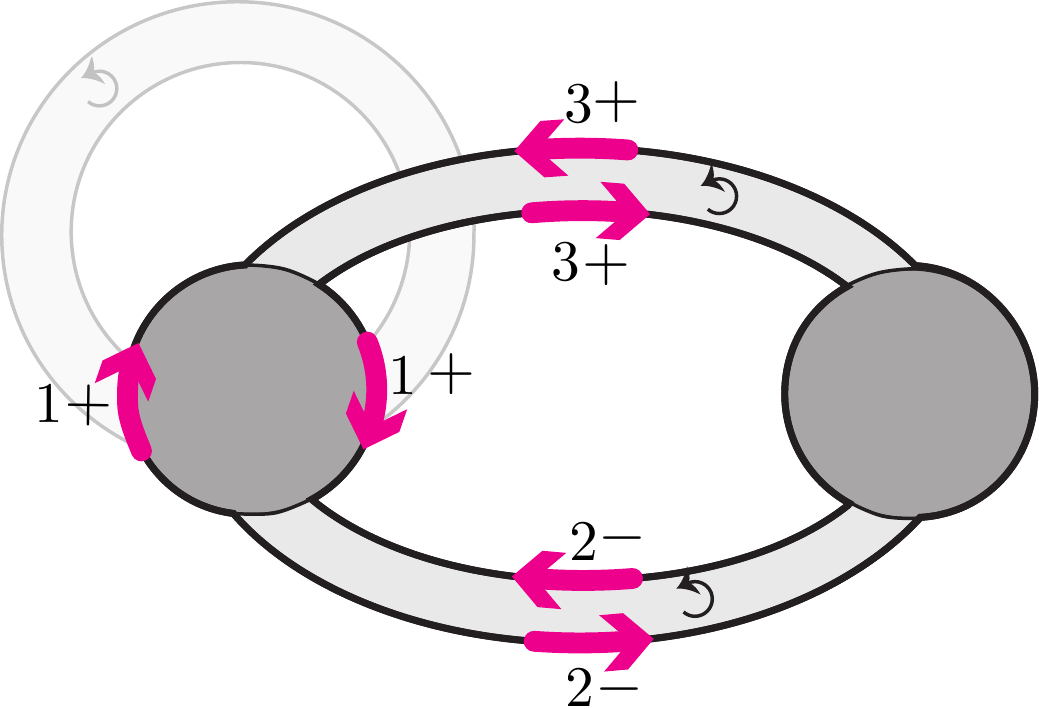}   \\
 Step 1.  & &Step 2. \ 
 \end{tabular}
 \end{center}

 \begin{center}
\begin{tabular}{ccc}
\includegraphics[width=5cm]{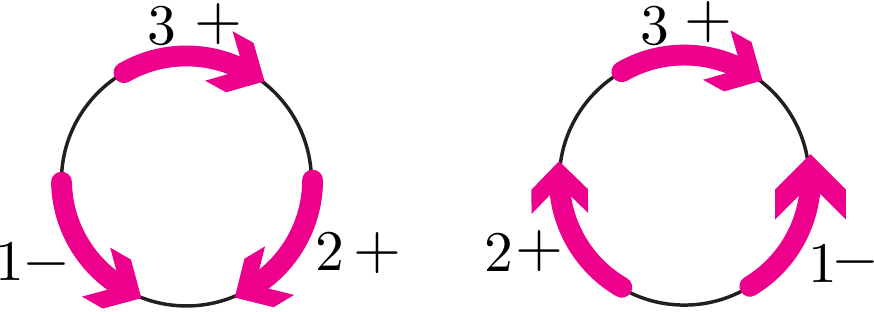} & \hspace{2cm} & \includegraphics[width=5cm]{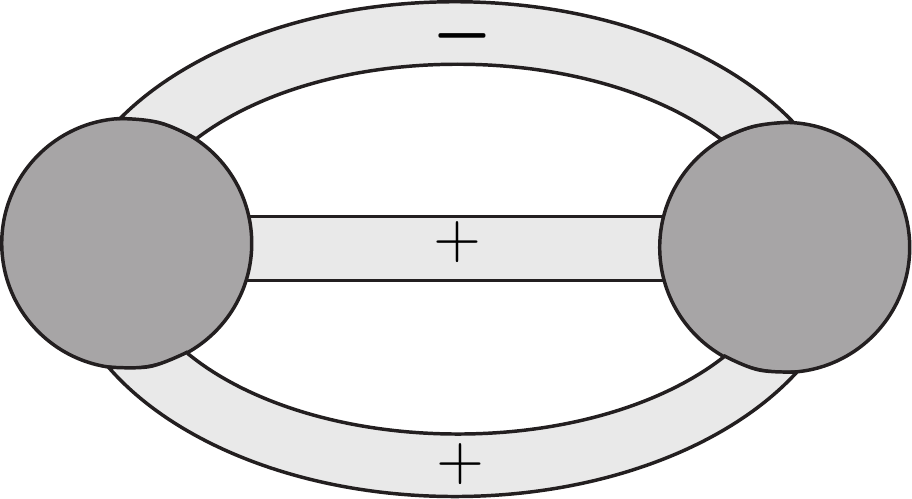}      \\
 Step 3.& & Steps  4.
\end{tabular}
\end{center}
\end{example}

\begin{example}\label{ex.dex2}
Again, the signed ribbon graph $G$ equipped with and arbitrary labelling and orientation of its edges is shown in Step 1. In this example $A=\{1,2\}$ and the marked spanning ribbon sub-graph $(V(G),A)$ is shown in Step 2. The resulting signed arrow presentation, shown in Step 3 and the partial dual   $G^{\{1,2\}}$ is shown in Step 4.  Note that in this example, $G$ and $G^A$ are equal as ribbon graphs, but not as signed ribbon graphs.
 \begin{center}
\begin{tabular}{ccccc}
 \includegraphics[width=3cm]{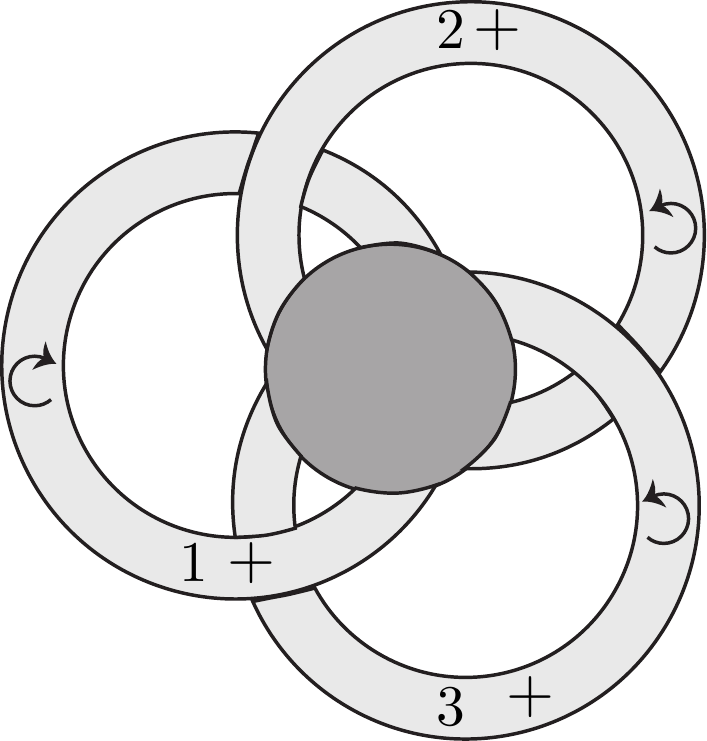} \hspace{1.5cm}& \includegraphics[width=3cm]{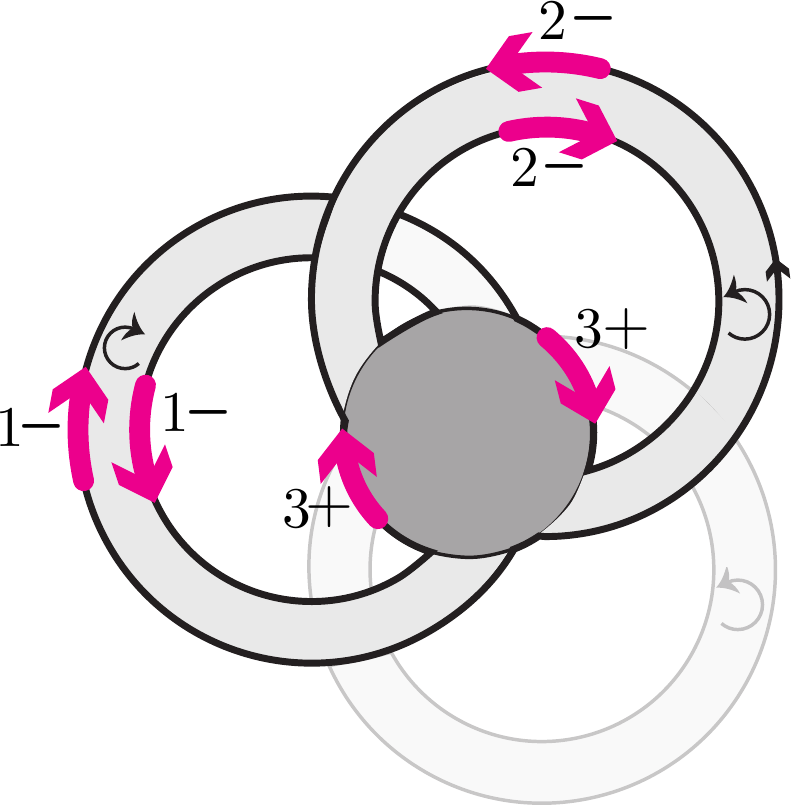}  \hspace{1cm} & \raisebox{3mm}{ \includegraphics[width=2cm]{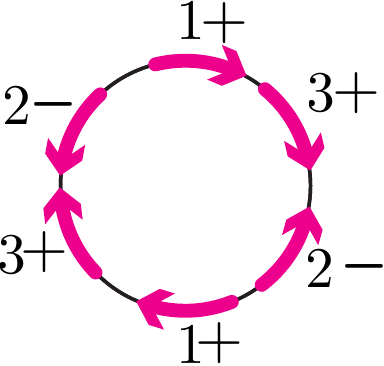} }\hspace{1cm} & \includegraphics[width=3cm]{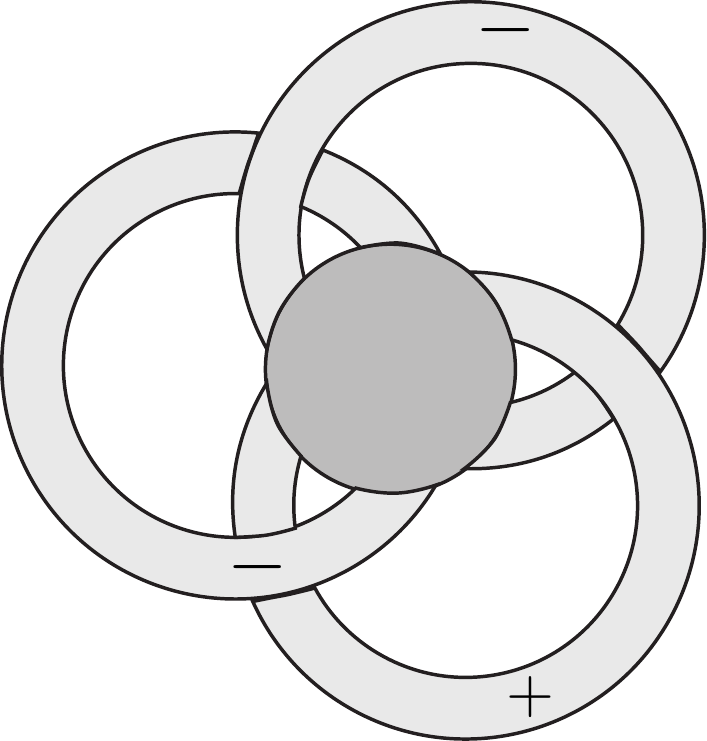}  \\
 Step 1. &   Step 2. & Steps 3 & Step 4. \\ 

\end{tabular}
 \end{center}
\end{example}

In both of these examples $G$ and $G^A$ have the same number of vertices. In general this will not be the case (for example the partial dual of a 2-cycle taken with respect to one edge is the non-planar, one vertex, two edge orientable  ribbon graph). Also notice that  $G$ and $G^A$ can have different genera. However, $G$ and $G^A$ will always have the same number of edges and the same number of connected components. In addition, in \cite{Ch1}, Chmutov observed that for all $A\subseteq E(G)$,    $G$ is orientable if and only if $G^A$ is orientable.

Observe that the definition of partial duality gives rise to a natural bijection between the edge set  $E(G)$ and $E(G^A)$. If $e$ is an edge of $G$, I will denote the corresponding edge in $G^A$ by $e^A$.

\begin{remark}
The dual $G^*$ of $G$ is formed in the following way: regarding $G$ as a punctured surface, fill in the punctures with disks and delete the original vertex set. The  resulting ribbon graph is $G^*$.  Chmutov observed in \cite{Ch1} that  $G^{E(G)}$ is the usual dual ribbon graph $G^*$ with all of the edge weights reversed. \end{remark}

\subsection{A geometric description of partial duals}\label{ss.gpd}

I will now provide a geometric description of the partial dual of a ribbon graph locally in the neighbourhood an edge. This geometric description will be especially convenient when we consider the homfly polynomial later.

Let $e$ be an edge of a signed ribbon graph $G$ and let $\varepsilon$ denote the sign of this edge.  We would like to know what the corresponding edge $e^A$ of the partial dual $G^A$ will look like.  There are two cases to consider: when $e\notin A$
 and when $e\in A$. We will deal with the easier case, $e\notin A$, first.
 
 \medskip
 
 Suppose that $e\notin A$. By untwisting the edge if necessary, we may assume that  the edge $e$ of the ribbon graph looks like
 \raisebox{-4mm}{\includegraphics[height=10mm]{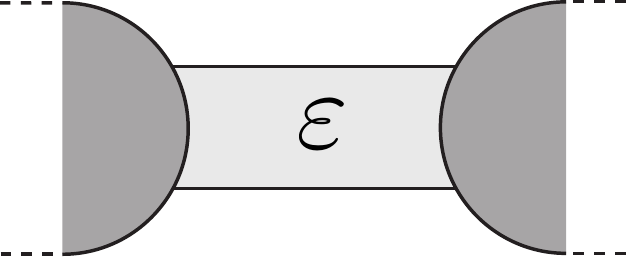}}. 
By the definition of partial duality, it follows that the edge $e^A$ is given by the arrow presentation  \raisebox{-4mm}{\includegraphics[height=10mm]{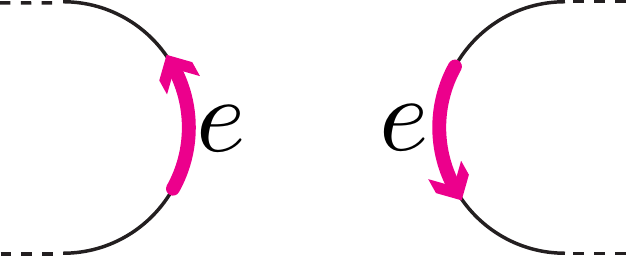}}, 
where the sign of the label $e$ in the arrow presentation is $\varepsilon$. Thus the edge $e^A$ of $G^A$ that corresponds to $e$ looks locally like \raisebox{-4mm}{\includegraphics[height=10mm]{local1}}.  That is, if $e\notin A$, we can assume that $G$ and $G^A$ are unchanged in a neighbourhood of the edge $e$.  This case can be summarized by the following table:

\begin{center}
\begin{tabular}{|c|c|}
\hline
 \raisebox{0mm}{\includegraphics[width=4cm]{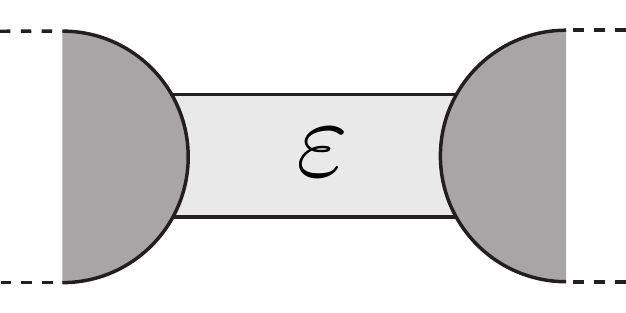}}  &  \raisebox{0mm}{\includegraphics[width=4cm]{r2}}  \\
 \hline
  $e \in E(G)$  &       $e^A\in E(G^A)$ when $e\notin A$ \\ 
\hline
\end{tabular}\;.
 \end{center}

\medskip

The second case, when $e\in A$, is a little more involved. Again, we may assume that  the edge $e$ of $G$ looks like  \raisebox{-4mm}{\includegraphics[height=10mm]{local1}}. By the definition of partial duality, it follows that the edge $e^A$, in the same neighbourhood, is given by the arrow presentation 
 \[\raisebox{-4mm}{\includegraphics[height=15mm]{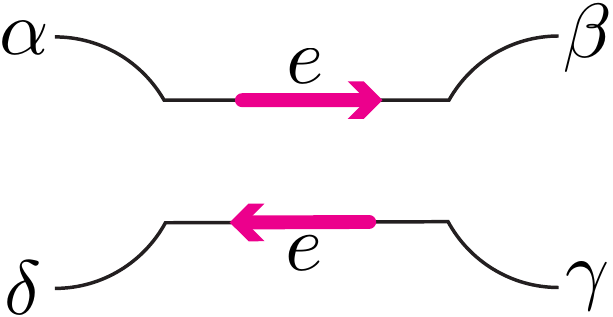}},\] 
where the sign of the label $e$ in the arrow presentation is $-\varepsilon$. Let $\alpha$, $\beta$, $\gamma$ and $\delta$ be the the points on the arcs of the arrow presentation shown in the figure. Then one of two things can happen: either $\alpha$, $\beta$, $\gamma$ and $\delta$ belong to the same cycle of the arrow presentation, or they do not. We will deal with each of these cases separately. 

\noindent \textit{Subcase 1.} If $\alpha$, $\beta$, $\gamma$ and $\delta$ all belong to the same cycle of the arrow presentation, then they  must appear in the cyclic order $(\alpha\, \beta\, \gamma\, \delta)$ or $(\alpha\,\gamma\, \delta\, \beta )$ with respect to some orientation of the cycle.  In either case we may assume that in our  drawing of $G^A$in the neighbourhood of  $e^A$, the single vertex incident with $e^A$ ``fills the gap'' left by the edge:
\[\includegraphics[width=4cm]{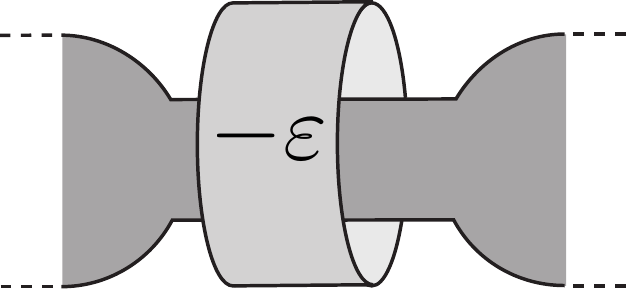}\;.\]

\noindent \textit{Subcase 2.} If $\alpha$, $\beta$, $\gamma$ and $\delta$ all belong different cycles of the arrow presentation, then  $\alpha$ and  $\beta$ lie on  one cycle, and  $\gamma$ and $\delta$ lie on  another cycle. Geometrically, this means that we can assume that  a neighbourhood of  the edge $e^A$ looks like 
\[ \includegraphics[height=3cm]{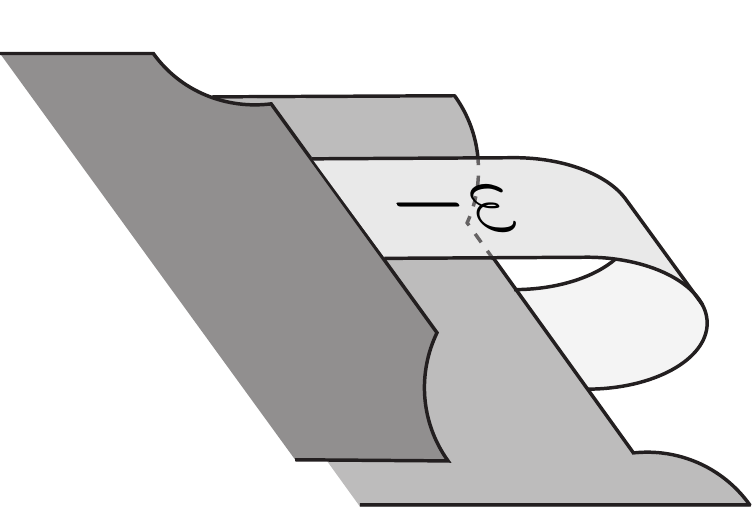}\]
(or a reflection in the vertical or in the plane on which it is drawn).
In this figure the edge is incident with two distinct vertices with the darker coloured vertex sitting above the lighter coloured vertex. Observe that the figure above can be deformed so as   to flatten out the edge:
 \[\raisebox{15mm}{\includegraphics[width=18mm]{arrow}}\includegraphics[height=3cm]{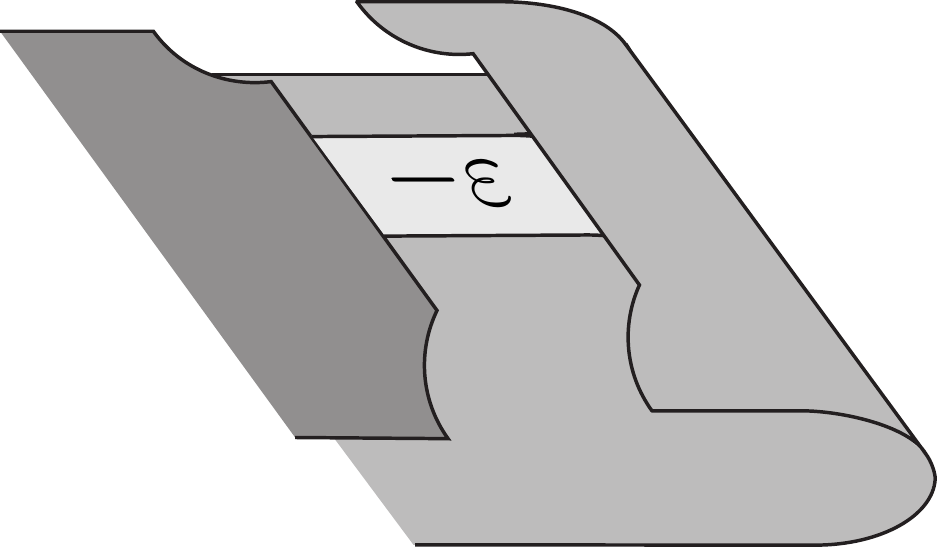} \raisebox{15mm}{\includegraphics[width=18mm]{arrow}} \includegraphics[height=3cm]{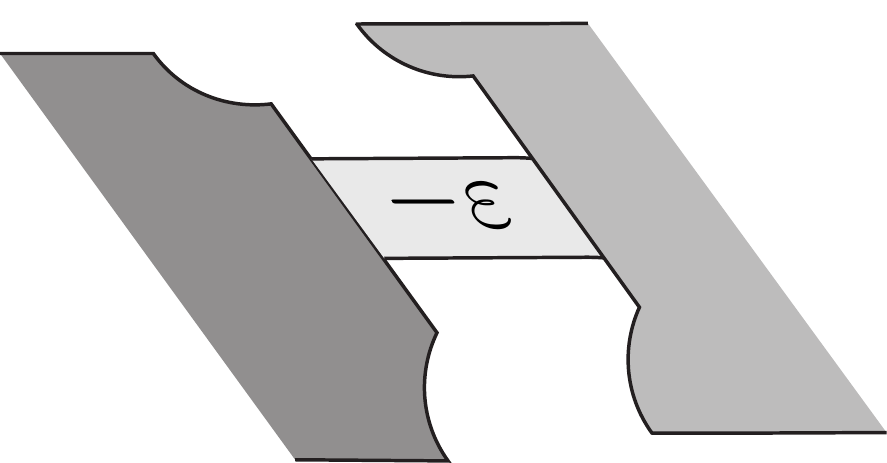}.\]
The figure on the left is obtained by ``straightening out the edge'' and the figure on the right is obtained from the figure on the left by taking a smaller neighbourhood of the edge.

This completes the  analysis of the  case when  $e\in A$. This analysis  is summarized by the following table.
  
\begin{center}
\begin{tabular}{|c|c|}
\hline
 \raisebox{-2mm}{\includegraphics[width=4cm]{r2} } & \includegraphics[width=4cm]{local2} \hspace{1cm} \raisebox{7mm}{or}     \hspace{1cm} \raisebox{-8mm}{ \includegraphics[height=3.2cm]{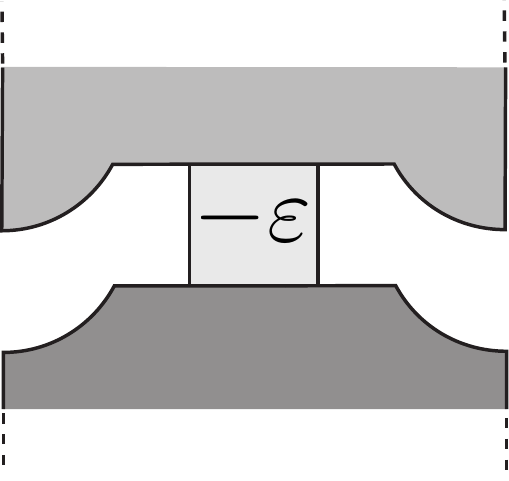}} \\
 \hline
  $e\in E(G)$  &        $e^A\in E(G)$ when $e\in A$ \\ 
\hline
\end{tabular}.
 \end{center}

\section{Polynomials}\label{s.poly}
\subsection{The signed ribbon graph polynomial}

I will begin by fixing some notation. Let $G$ be a signed ribbon graph with vertex set $V(G)$ and edge set $E(G)$. Let 
$v(G)=|V(G)|$, $e(G)=|E(G)|$, $k(G)$ be the number of connected components of $G$, $\partial (G)$ be the number of boundary components of $G$, $r(G)=v(G)-k(G)$ and $n(G)=e(G)-r(G)$. In addition, let $e_+(G)$ denote the number of positively signed edges of $G$, and   $e_-(G)$ denote the number of negatively signed edges of $G$. Finally, a {\em state} of a signed ribbon graph is a signed spanning ribbon sub-graph. (So a state of $G$ is a signed ribbon graph found from $G$ by removing a subset of edges.)  Let $\F (G)$ denote the set of the $2^{E(G)}$ states of $G$.

\medskip

The {\em signed ribbon graph  polynomial} was introduced by Chmutov and Pak in \cite{CP}.  Along with its generalizations it has appeared in several places in the literature (for example \cite{CV,Ch1,HM,LM,Mo1,Mo2}). It is defined  by the state sum
\begin{equation}
\label{eq.brpoly}
 R_s(G\; ; x,y,z) = \sum_{F\in \F (G)}  x^{r(G)-r(F)+s(F)} y^{n(F)-s(F)}z^{k(F)-\partial(F)+n(F)}
\end{equation}
where 
\[  s(F)  = \frac{1}{2}(e_-(F) - e_-(G-F)).   \]

The  signed ribbon graph  polynomial is an element of $\mathbb{Z} [ x^{\pm \frac{1}{2}},  y^{\pm \frac{1}{2}}, z^{\pm 1} ]$.

\begin{example}\label{ex.wbr}
The signed ribbon graph $G$ from example~\ref{ex.dex1} has the signed ribbon graph  polynomial
\[   R_s(G\; ; x,y,z) =  x^{\frac{1}{2}}y^{\frac{3}{2}}z^2+  x^{\frac{1}{2}}y^{\frac{3}{2}} + 3  x^{\frac{1}{2}}y^{\frac{1}{2}} +  x^{-\frac{1}{2}}y^{\frac{3}{2}} + x^{-\frac{1}{2}}y^{\frac{1}{2}}+ x^{\frac{1}{2}}y^{-\frac{1}{2}}, 
\]
and for the signed ribbon graph $G^A$ from the same example
\[   R_s(G^A\; ; x,y,z) =    x^{\frac{1}{2}}y^{\frac{3}{2}} + 3  x^{\frac{1}{2}}y^{\frac{1}{2}} + x^{\frac{1}{2}}y^{-\frac{1}{2}} +  x^{-\frac{1}{2}}y^{\frac{3}{2}} +2 x^{-\frac{1}{2}}y^{\frac{1}{2}}.
\]
\end{example}

I can now write down Chmutov's duality theorem.
\begin{theorem}[Chmutov \cite{Ch1}]\label{th.main}
If $G$ is a signed ribbon graph and $G^A$ is a partial dual of $G$, then when $xyz^2=1$,
\begin{equation}\label{th.ch}
 (yz)^{v(G)} R_s\left( G ; x,y,z\right) =   (yz)^{v(G^A)} R_s\left( G^A ; x,y,z\right). 
\end{equation}
\end{theorem}
As I have mentioned previously, the aim of this paper is to provide a new and simple proof of this theorem for orientable ribbon graphs through the use of basic knot theory.

Notice that example~\ref{ex.wbr} verifies this theorem.

\subsection{The homfly polynomial}
The homfly polynomial \cite{homfly,PT} of a link in ${\bf S}^3$  (or $\mathbb{R}^3$)  can be  defined recursively by the relations
\begin{equation}\label{eq.h1}
 X \, P\left( L_+ \right)-X^{-1} \,  P\left( L_- \right) =Y\, P\left( L_0 \right) \end{equation}
 and
 \begin{equation}
 \label{eq.h2}  P\left( \mathcal{O}^k \right)=\left(\frac{X-X^{-1}}{Y}\right)^{k-1},
\end{equation}
where $\mathcal{O}^k$ is a $k$ component unlink diagram ({\em i.e.} the $k$ component link with no crossings), and $L_+$, $L_-$ and $L_0$ are link diagrams which are identical except in a single region where they differ as indicated:

\begin{center}
\begin{tabular}{ccccc}
 \includegraphics[width=1.5cm]{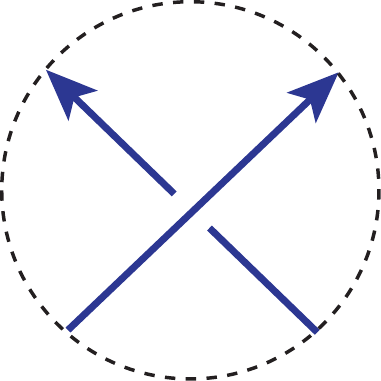} & \hspace{1cm} & \includegraphics[width=1.5cm]{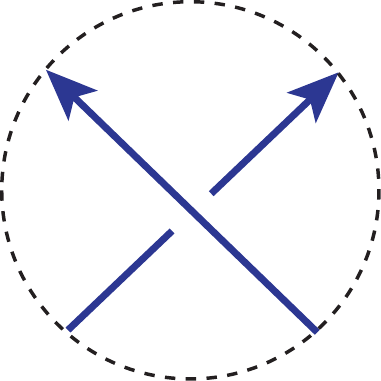} & \hspace{1cm} &   \includegraphics[width=1.5cm]{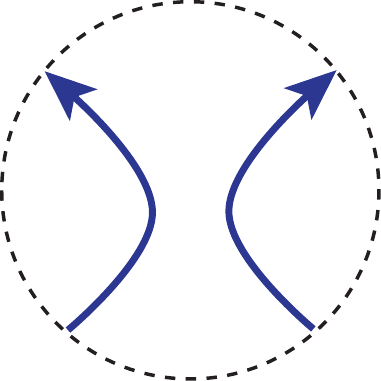}  \\ & &&& \\
$L_+$ &  & $L_-$ & & $L_0$
 \end{tabular}.
\end{center}

In ${\bf S}^3$, the relations \ref{eq.h1} and \ref{eq.h2} define a polynomial since the skein relation \ref{eq.h1} will reduce any link diagram on ${\bf S}^2$ to a $\mathbb{Z}[X^{\pm 1}, Y^{\pm 1}]$ linear combination of unlink diagrams. Equation~\ref{eq.h2} can then be used to obtain a polynomial. 

For link diagrams on an arbitrary orientable surface, however, the skein relation \ref{eq.h1}  will not necessarily reduce a link diagram to a linear combination of unlink diagrams, so equations \ref{eq.h1} and \ref{eq.h2} are not enough to  define the homfly polynomial. A little more work is required to define the homfly polynomial of a link diagram on an arbitrary orientable surface. It was shown in \cite{Lie} that the skein relation \ref{eq.h1} will reduce a link to a linear combination of {\em descending links} (a definition of descending links will follow shortly). A homfly polynomial of a link diagram on a surface can then be defined by specifying its values on descending links. Here I will set 
\begin{equation}\label{eq.h3}
P\left( \mathcal{D} \right)=\left(\frac{X-X^{-1}}{Y}\right)^{k-1},
\end{equation}
where $\mathcal{D}$ is a $k$ component descending link, and define the {\em homfly polynomial} to be the unique polynomial defined by equations \ref{eq.h1} and \ref{eq.h3}. (More general multivariate homfly polynomials can be defined  by choosing a basis for the homfly skein that depends on the conjugacy class of the descending links in the fundamental group of the surface, see \cite{Lie} for details. Here, however, this extra generality is not needed.)

I will now give a definition of a descending link. The following concept of a product is needed for the definition of a descending link. Let $\Sigma$ be an orientable surface.
There is a natural {\em product} of links in $\Sigma\times I$ given by reparameterizing the two copies of  $\Sigma\times I$  and stacking them:
\[
 (\Sigma \times I) \times  (\Sigma \times I) \cong (\Sigma \times [1/2,1]) \times  (\Sigma \times [0,1/2])
\rightarrow (\Sigma \times I) .
\]
Also denote the projections from $\Sigma \times I$ to $\Sigma $ and to $I$ by $p_\Sigma $ and $p_I$ respectively.  The value $p_I(x)$ is called the {\em height} of $x$. 
\begin{definition}

(1) A knot $K \subset \Sigma  \times I$ is {\em descending} if
it is isotopic to a knot $K^{\prime} \subset \Sigma  \times I$ with the property that
 there is a choice of basepoint $a$ on $K^{\prime}$ such that if we travel along $K^{\prime}$ in the direction of the orientation from the basepoint the height of $K^{\prime}$ decreases until we reach a point $a^{\prime}$ with $p_\Sigma (a)=p_\Sigma (a^{\prime})$ from which $K^{\prime}$ leads back to $a$ by increasing the height and keeping the projection onto $F$ constant.

(2) A link $L \subset \Sigma  \times I$ is said to be {\em descending} if it is isotopic to a product of descending knots.
\end{definition}

The following example will be important later.
\begin{example}\label{e.desc}
Any link diagram on an orientable surface that has no crossings  is a diagram of a descending link.
\end{example}

\subsection{The homfly and the ribbon graph  polynomial}\label{ss.hbr}
In \cite{Mo1} I described a relation between the homfly polynomial of a certain class of links in thickened surfaces and the ribbon graph  polynomial. This relation generalized earlier results of  Jaeger \cite{Ja} and Traldi \cite{Tr} which relate the homfly polynomial of a link in    ${\bf S}^3$ with the Tutte polynomial of a planar graph.
I will use the connection between the homfly and ribbon graph  polynomials to prove Chmutov's duality theorem. The relevant property from \cite{Mo1} is as follows: 
given an orientable  signed ribbon graph $G$, construct a link diagram $\calL_G$ on $G$ by associating the following configurations at each signed edge of $G$  
\begin{center}
\begin{tabular}{ccc}
 \includegraphics[width=4cm]{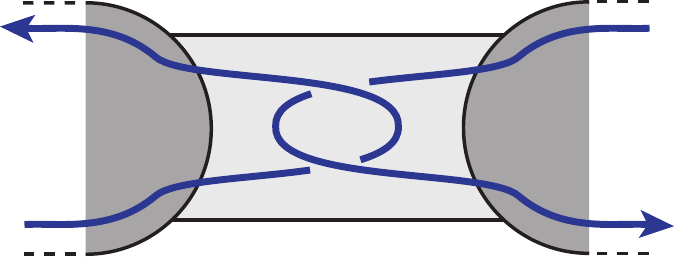} & \hspace{2cm} &   \includegraphics[width=4cm]{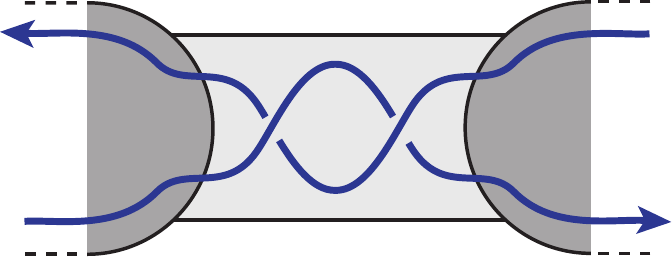}  \\ & & \\
 $\calL_G$ at a $+$ edge & & $\calL_G$ at a $-$ edge
 \end{tabular}
\end{center}
and connecting the configurations by following the boundary of the vertices.   This gives a diagram of a link in the thickened surface $G \times I$ where $I=[0,1]$ is the unit interval. 
\begin{example}
 If $G$ is the ribbon graph from example~\ref{ex.dex1}, then $\calL_{G}$ is the link diagram
\[ \includegraphics[width=6cm]{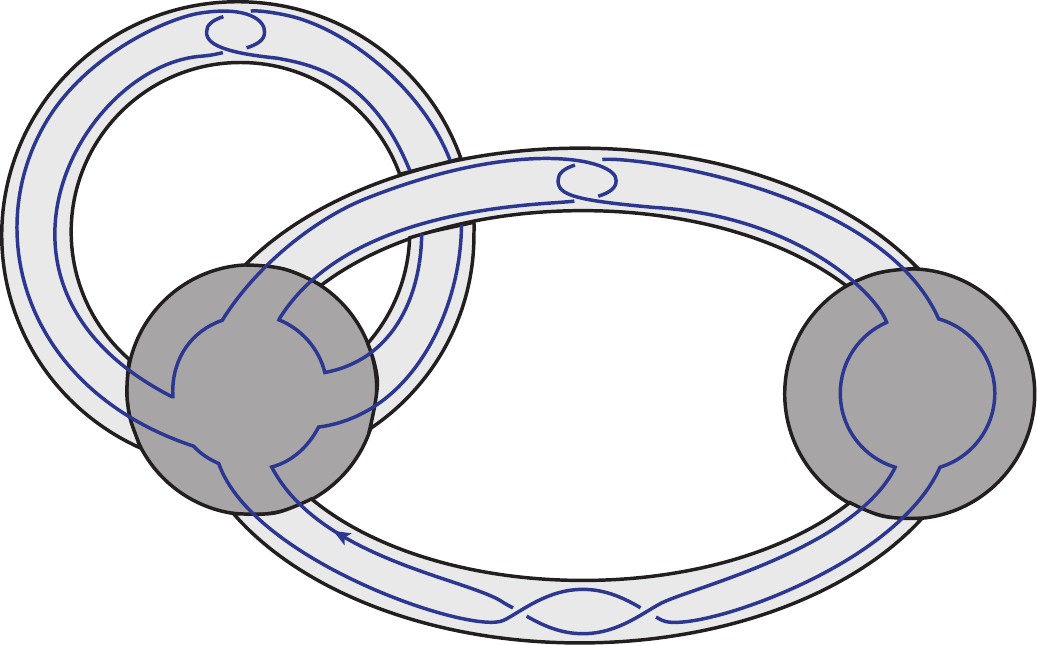}. \]
 \end{example}

It was shown in \cite{Mo1} that if $G$ is an orientable signed ribbon graph, then the homfly polynomial of the link $\calL_G$ is an evaluation of the multivariate  ribbon graph  polynomial. In fact, 
 Theorem~4.3 of \cite{Mo1} gives
\begin{equation}\label{eq.brhom}
P\left( \calL_G ; X,Y   \right) = 
\left(  \frac{Y}{X-X^{-1}}   \right)\left(  \frac{Y}{X}   \right)^{e_- (G)}\left(  \frac{1}{X^2}   \right)^{e_+ (G)}
\sum_{F\in \F (G)} \left(  \frac{X-X^{-1}}{Y}\right)^{\partial (F)} \prod_{e\in F} w_e,
\end{equation}
where 
 \[ 
 w_e = \left\{ \begin{array}{ll}   XY & \quad \text{ if }  e \text{ of positive weight,}  \\   \frac{1}{XY} & \quad \text{ if }  e \text{ of negative weight.}  \end{array}   \right. 
\]
I will use this identity to reduce Chmutov's duality theorem (Theorem~\ref{th.main}) to a simple knot theoretic problem.

\section{A proof of Chmutov's duality theorem}\label{s.proof}
In this section all of our ribbon graphs $G$ will be orientable.
\subsection{A knot theoretic reformulation}\label{ss.kt}

Expanding the rank and nullity in equation~\ref{eq.brpoly} and collecting terms gives 
\[  R_s(G;\; x,y,z) = x^{-k(G)}  (yz)^{-v(G)} \sum_{F\in \F (G)}  (xyz^2)^{k(F)} (yz)^{e(F)}z^{-\partial(F)} (xy^{-1})^{s(F)} . \] 
Making the substitutions $a=xyz^2$, $b=zy$ and $c=z^{-1}$ then gives
\[  R_s\left(G; \;\frac{ac}{b},bc, \frac{1}{c}\right) =  \left( \frac{b}{ac}\right)^{k(G)}  \left( \frac{1}{b}\right)^{v(G)} \sum_{F\in \F (G)}  a^{k(F)} b^{e(F)}   \left( \frac{a}{b^2}\right)^{s(F)}  c^{\partial(F)}. \]
We now turn our attention to rewriting the sign function $s(F)$ in this expression. The term $e_- (G-F)$ used in the definition of $s(F)$ can be expressed as 
$ e_-(G_F) = e(G) - e_+(G)-e_-(F)$. 
Substituting this into the formula for $s(F)$ gives
\[ s(F)=\frac{1}{2}[ e_-(F)-e(G)+e_+(G)+e_-(F) ]=\frac{1}{2}[ 2e_-(F)-e_-(G) ], \]
where the second equality follows since $e(G)-e_+(G)=e_-(G)$.
Thus we have 
\[ b^{e(F)} \left(\frac{a}{b^2}\right)^{s(F)} =
 b^{e(F)} \left(\frac{a}{b^2}\right)^{e_-(F)-\frac{1}{2}e_-(G)}  =   
 \left( \frac{b}{a^{1/2}}\right)^{e_-(G)}  \left[a^{e_-(F)} b^{e(F)-2e_-(F)} \right]
 =  \left( \frac{b}{a^{1/2}}\right)^{e_-(G)}  \prod_{e\in F} \omega_e,\]
 where
 \[ 
 \omega_e = \left\{ \begin{array}{ll}   b & \quad \text{ if }  e \text{ of positive weight,}  \\   \frac{a}{b} & \quad \text{ if }  e \text{ of negative weight.}  \end{array}   \right. 
\]
We can now write the signed ribbon graph  polynomial as a Potts model type state sum:
\begin{equation}\label{eq.brpotts}
  R_s\left(G; \;\frac{ac}{b},bc, \frac{1}{c}\right) = \left( \frac{b}{ac}\right)^{k(G)}  \left( \frac{1}{b}\right)^{v(G)} \left( \frac{b}{\sqrt{a}}\right)^{e_-(G)} \sum_{F\in \F (G)}  a^{k(F)}  c^{\partial(F)} \prod_{e\in F} \omega_e. 
\end{equation}
Now setting $a=1$ in equation~\ref{eq.brpotts}, and $X=\sqrt{bc+1}$ and $Y=\frac{b}{\sqrt{bc+1}}$ into equation~\ref{eq.brhom}, the sums on the right hand side of the two expressions equate and we can write
\begin{multline*}
P\left( \calL_G ; \sqrt{bc+1},\frac{b}{\sqrt{bc+1}}   \right)  \\ = 
\frac{1}{c} \left(   \frac{b}{bc+1} \right)^{e_-(G)} \left(   \frac{1}{bc+1} \right)^{e_+(G)} b^{e_-(G)}\left(   \frac{c}{b} \right)^{k(G)} b^{v(G)} \left(   \frac{1}{b} \right)^{e_-(G)}
 R_s\left(G; \;\frac{c}{b},bc, \frac{1}{c}\right)   \\
 =  \frac{1}{c}  \left(   \frac{1}{bc+1} \right)^{e(G)} \left(   \frac{c}{b} \right)^{k(G)} b^{v(G)} R_s\left(G; \;\frac{c}{b},bc, \frac{1}{c}\right).
\end{multline*}
Finally, recovering the original variables $x$, $y$ and $z$ using  $a=xyz^2=1$, $b=zy$ and $c=z^{-1}$ and simplifying, gives the identity
\begin{equation}\label{eq.pr}
P\left( \calL_G ; \sqrt{y+1},\frac{yz}{\sqrt{y+1}}   \right) =
(y+1)^{-e(G)}
x^{k(G)}y^{v(G)}z^{v(G)+1}  R_s(G;\; x,y,z) ,
\end{equation}
where $xyz^2=1$.

By substituting equation~\ref{eq.pr} in to the left and right hand sides of equation~\ref{th.ch}, we obtain the following reformulation of theorem~\ref{th.main}.
\begin{lemma}\label{l.knot}
Theorem~\ref{th.main} holds if and only if 
\begin{equation}\label{eq.lem}
P\left( \calL_G ;X,Y  \right) =
P\left( \calL_{G^A} ; X,Y  \right),
\end{equation}
for all orientable signed ribbon graphs $G$ and for all $A\subseteq E(G)$.
\end{lemma}
I give a straightforward proof of this lemma, and therefore of Chmutov's duality theorem, in the following subsection.

\subsection{A proof of the theorem}\label{ss.proof}

In this final subsection I prove that equation~\ref{eq.lem} does indeed hold, and thus, by lemma~\ref{l.knot}, theorem~\ref{th.main} also holds. I prove equation~\ref{eq.lem} by considering the contributions of the links $\calL_G$ and $\calL_{G^A}$  at an edge $e$ of the ribbon graphs $G$ and $G^A$ (recall that there is a bijection between the edges of $G$ and the edges of $G^A$) to the homfly polynomial.
\begin{lemma}
Let $G$ be an orientable  signed ribbon graph. Then
\begin{equation*}
P\left( \calL_{G} ;X,Y  \right) =
P\left( \calL_{G^A} ; X,Y  \right),
\end{equation*}
 for all $A\subseteq E(G)$.
\end{lemma}
\begin{proof}
Let $e$ be an edge of $G$. First suppose that $e$ is of positive weight. Then the contribution at $e$ to $P\left( \calL_{G} ;X,Y  \right)$ is calculated as follows
\[
\begin{array}{cll}
 \includegraphics[width=3cm]{Lplus} &\raisebox{5mm}{$ =  \frac{1}{X^2}$} \includegraphics[width=3cm]{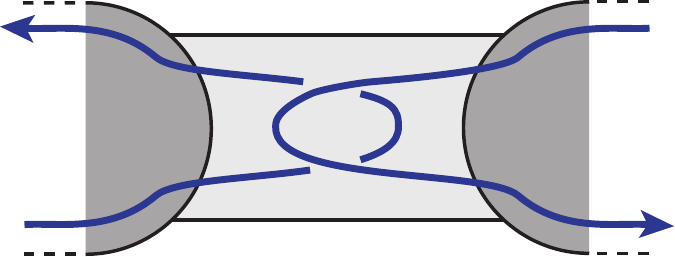} & \raisebox{5mm}{$  +\frac{Y}{X}$} \includegraphics[width=3cm]{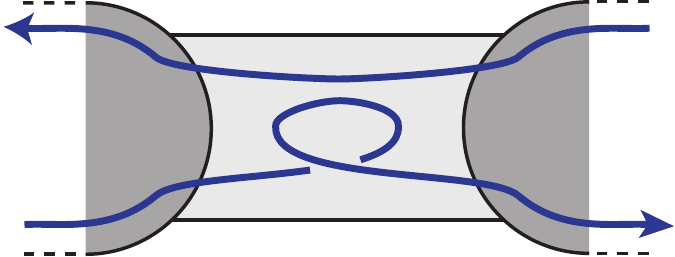}
  \\ && \\
 &\raisebox{5mm}{ $=     \frac{1}{X^2} $}\includegraphics[width=3cm]{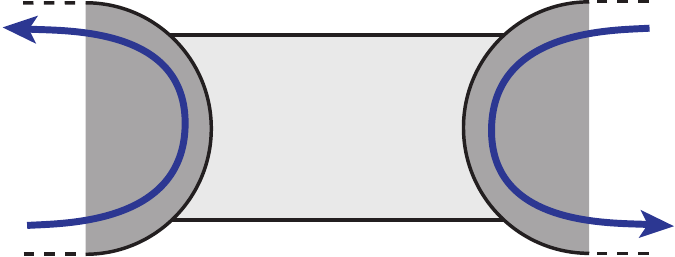} &   \raisebox{5mm}{$+\frac{Y}{X}$} \includegraphics[width=3cm]{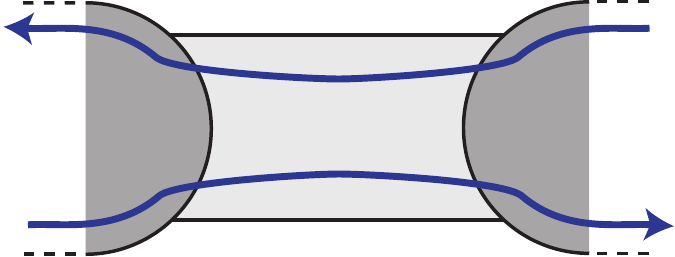}
 \end{array},
\]
where the first equality follows from \ref{eq.h1} and the second follows by isotopy.

Now consider the contribution at the corresponding edge $e$ in $G^A$. If $e \notin A$, then, as described in the first table in Subsection~\ref{ss.gpd}, locally the edge $e$ in $G^A$ is the same as the edge $e$ in $G$. This means that locally at $e$ the links $\calL_{G}$ and $\calL_{G^A}$ are identical and so their contributions to the homfly polynomials are identical.

If $e\in A$, then at the edge $e$, $G^A$ and $G$ differ as described in the second table in Subsection~\ref{ss.gpd}. The homfly polynomial calculation for $\calL_{G^A}$ at $e$ is then
\[
\begin{array}{cll}
 \includegraphics[width=3cm]{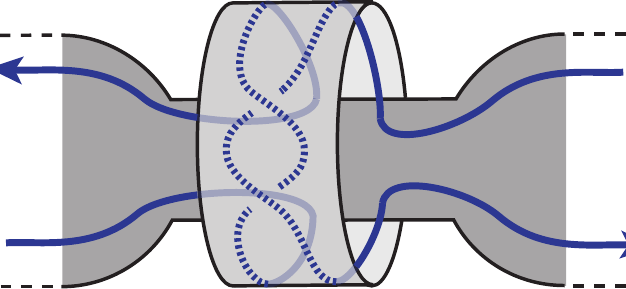} &\raisebox{5mm}{$ =  \frac{1}{X^2}$} \includegraphics[width=3cm]{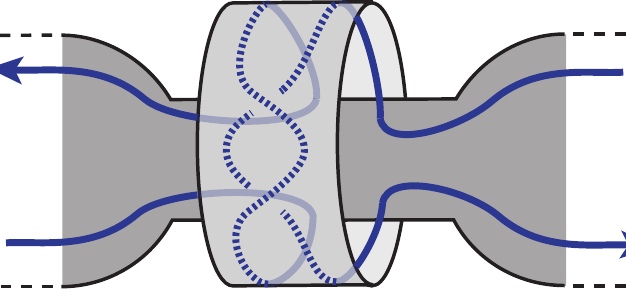} & \raisebox{5mm}{$  +\frac{Y}{X}$} \includegraphics[width=3cm]{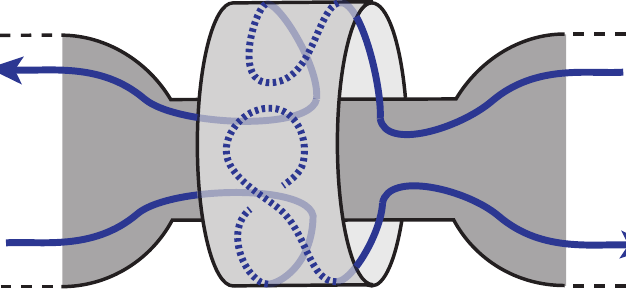} \\ && \\
 &\raisebox{5mm}{ $=     \frac{1}{X^2} $}\includegraphics[width=3cm]{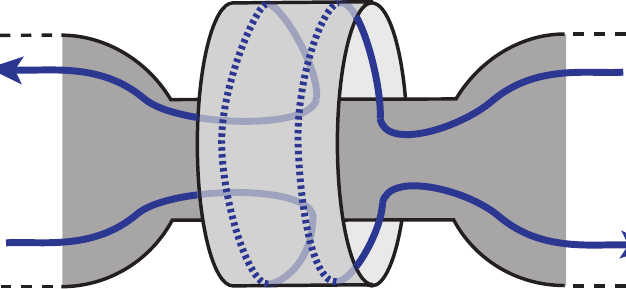} &   \raisebox{5mm}{$+\frac{Y}{X}$} \includegraphics[width=3cm]{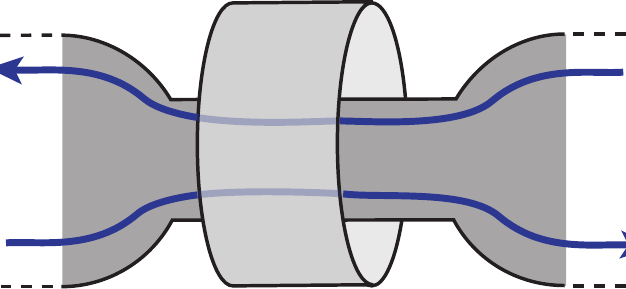}
 \end{array},
\]
or 
\[
\begin{array}{cllll}
 \includegraphics[height=3cm=3cm]{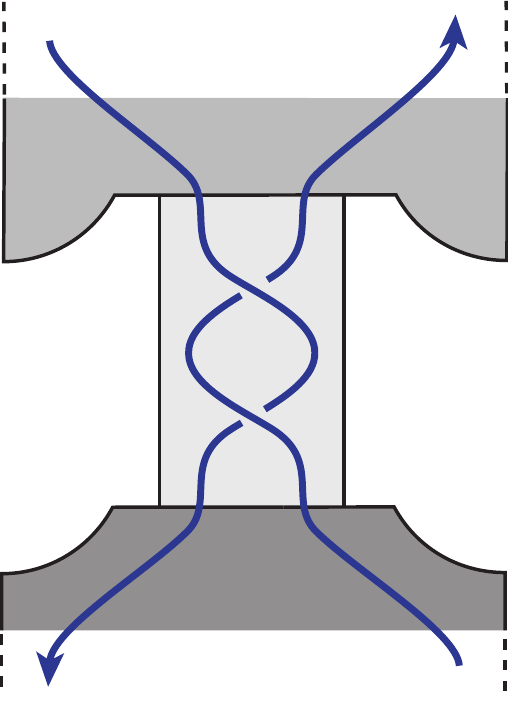} &\raisebox{15mm}{$ =  \frac{1}{X^2}$} \includegraphics[height=3cm]{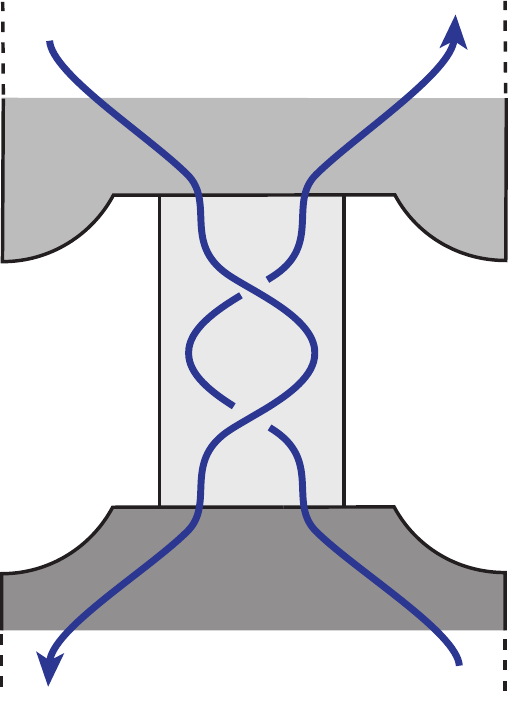} & \raisebox{15mm}{$  +\frac{Y}{X}$} \includegraphics[height=3cm]{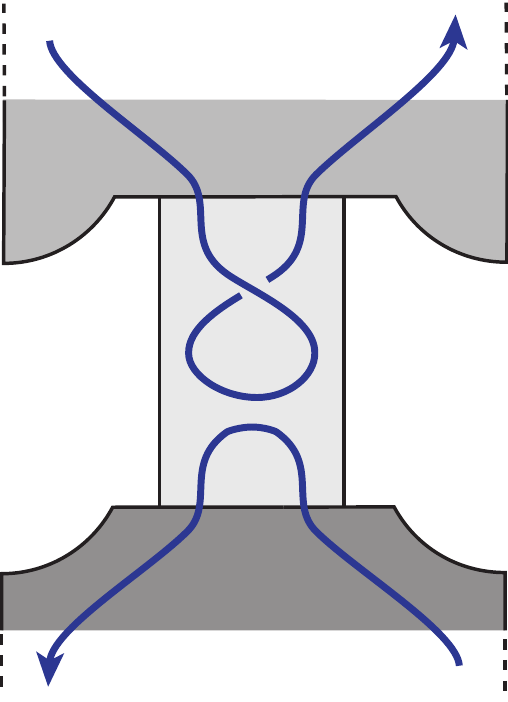} 
 &\raisebox{15mm}{ $=     \frac{1}{X^2} $}\includegraphics[height=3cm]{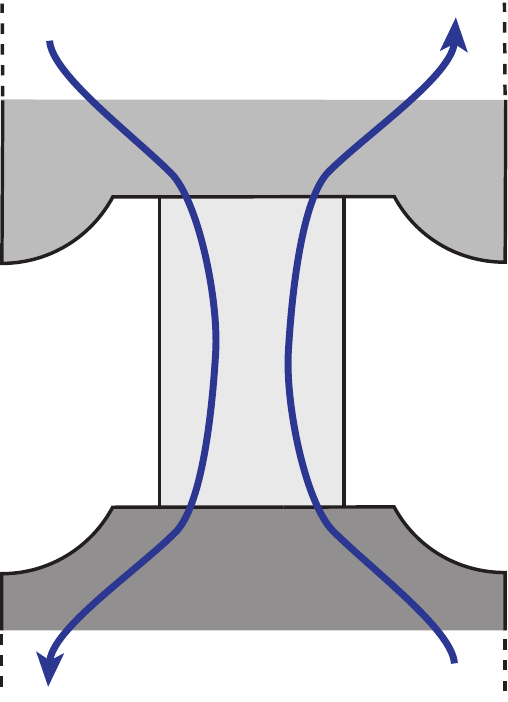} &   \raisebox{15mm}{$+\frac{Y}{X}$} \includegraphics[height=3cm]{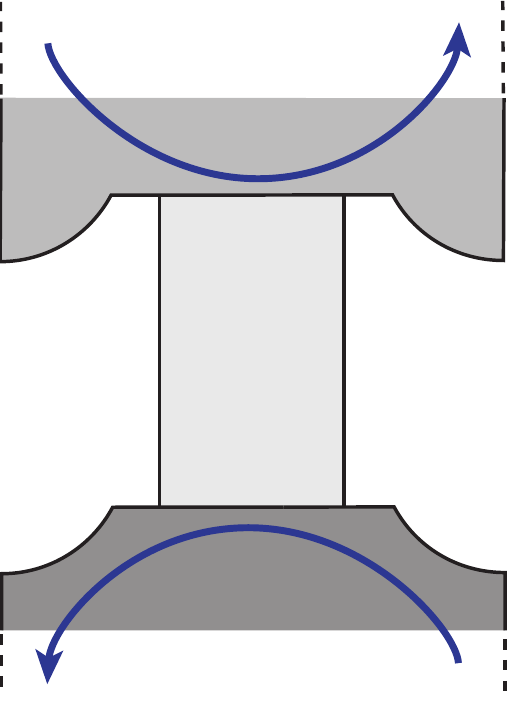}
 \end{array}.
\]


In the above example, the edge $e$ of $G$ was positive. A similar calculation can be done when the edge is negative.

Resolving every crossing of  $\calL_G$ and $\calL_{G^A}$ as indicated above gives two linear combinations of collections of cycles on the surfaces $G$ and $G^A$.   From the figures above, there is an obvious correspondence between the summands of the two linear combinations. Moreover,  the corresponding summands will  have the same number of cycles and the same coefficient. Finally, since the cycles in each summand do not cross, the cycles form a set of  descending links (by example~\ref{e.desc}), and we can calculate the homfly polynomials using \ref{eq.h3}. It then follows that $P\left( \calL_{G} ;X,Y  \right) =
P\left( \calL_{G^A} ; X,Y  \right)$ as required. 
\end{proof}

\begin{remark}
One can also use the knot theoretic approach above to prove Chmutov's change of sign formula proposition~2.5 of \cite{Ch1} along the surface $xyz^2=1$.

Also note that  Ellis-Monaghan and I. Sarmiento's duality relation for the ribbon graph  polynomial from \cite{ES} and \cite{Mo1} is a consequence of the above fact and  Chmutov's duality relation. See \cite{Ch1} Section~4.1 for details. \end{remark}



\end{document}